\documentclass[a4paper, 11pt]{article}
\usepackage[ansinew]{inputenc}
\usepackage{amssymb, amsmath, amsthm}
\usepackage{nicefrac, mathdots}
\usepackage{hyperref}

\def\typeout#1{\message{^^J}\message{#1}\message{^^J}}
\newif\ifSRCOK \SRCOKtrue
\newcount\PAGETOP
\newcount\LASTLINE
\global\PAGETOP=1
\global\LASTLINE=-1
\def\EJECT{\SRC\eject}
\def\WinEdt#1{\typeout{:#1}}
\gdef\MainFile{\jobname.tex}
\gdef\CurrentInput{\MainFile}
\newcount\INPSP
\global\INPSP=0
\def\SRC{\ifSRCOK%
  \ifnum\inputlineno>\LASTLINE%
    \ifnum\LASTLINE<0%
      \global\PAGETOP=\inputlineno%
    \fi%
    \global\LASTLINE=\inputlineno%
    \ifnum\INPSP=0%
      \ifnum\inputlineno>\PAGETOP%
        
      \fi%
    \else%
      
    \fi%
  \fi%
\fi}
\def\PUSH#1{%
\SRC%
\ifnum\INPSP=0 \global\let\INPSTACKA=\CurrentInput \else%
\ifnum\INPSP=1 \global\let\INPSTACKB=\CurrentInput \else%
\ifnum\INPSP=2 \global\let\INPSTACKC=\CurrentInput \else%
\ifnum\INPSP=3 \global\let\INPSTACKD=\CurrentInput \else%
\ifnum\INPSP=4 \global\let\INPSTACKE=\CurrentInput \else%
\ifnum\INPSP=5 \global\let\INPSTACKF=\CurrentInput \else%
               \global\let\INPSTACKX=\CurrentInput \fi\fi\fi\fi\fi\fi%
\gdef\CurrentInput{#1}%
\WinEdt{<+ \CurrentInput}%
\global\LASTLINE=0%
\ifSRCOK\fi%
\global\advance\INPSP by 1}
\def\POP{%
\ifnum\INPSP>0 \global\advance\INPSP by -1  \fi%
\ifnum\INPSP=0 \global\let\CurrentInput=\INPSTACKA \else%
\ifnum\INPSP=1 \global\let\CurrentInput=\INPSTACKB \else%
\ifnum\INPSP=2 \global\let\CurrentInput=\INPSTACKC \else%
\ifnum\INPSP=3 \global\let\CurrentInput=\INPSTACKD \else%
\ifnum\INPSP=4 \global\let\CurrentInput=\INPSTACKE \else%
\ifnum\INPSP=5 \global\let\CurrentInput=\INPSTACKF \else%
               \global\let\CurrentInput=\INPSTACKX \fi\fi\fi\fi\fi\fi%
\WinEdt{<-}%
\global\LASTLINE=\inputlineno%
\global\advance\LASTLINE by -1%
\SRC}
\def\INPUT#1{\relax}

\let\originalxxxeverypar\everypar
\newtoks\everypar
\originalxxxeverypar{\the\everypar\expandafter\SRC}
\everymath\expandafter{\the\everymath\expandafter\SRC}
\output\expandafter{\expandafter\SRCOKfalse\the\output}

\newif\ifSRCOK \SRCOKtrue
\newcount\PAGETOP
\newcount\LASTLINE
\global\PAGETOP=1
\global\LASTLINE=-1
\gdef\MainFile{\jobname.tex}
\gdef\CurrentInput{\MainFile}
\newcount\INPSP
\global\INPSP=0
\def\EJECT{\SRC\eject}
\def\WinEdt#1{\typeout{:#1}}
\def\SRC{\ifSRCOK%
  \ifnum\inputlineno>\LASTLINE%
    \ifnum\LASTLINE<0%
      \global\PAGETOP=\inputlineno%
    \fi%
    \global\LASTLINE=\inputlineno%
    \ifnum\INPSP=0%
      \ifnum\inputlineno>\PAGETOP%
      \fi%
    \else%
    \fi%
  \fi%
\fi}
\def\PUSH#1{%
\SRC%
\ifnum\INPSP=0 \global\let\INPSTACKA=\CurrentInput \else%
\ifnum\INPSP=1 \global\let\INPSTACKB=\CurrentInput \else%
\ifnum\INPSP=2 \global\let\INPSTACKC=\CurrentInput \else%
\ifnum\INPSP=3 \global\let\INPSTACKD=\CurrentInput \else%
\ifnum\INPSP=4 \global\let\INPSTACKE=\CurrentInput \else%
\ifnum\INPSP=5 \global\let\INPSTACKF=\CurrentInput \else%
               \global\let\INPSTACKX=\CurrentInput \fi\fi\fi\fi\fi\fi%
\gdef\CurrentInput{#1}%
\WinEdt{<+ \CurrentInput}%
\global\LASTLINE=0%
\ifSRCOK\fi%
\global\advance\INPSP by 1}
\def\POP{%
\ifnum\INPSP>0 \global\advance\INPSP by -1  \fi%
\ifnum\INPSP=0 \global\let\CurrentInput=\INPSTACKA \else%
\ifnum\INPSP=1 \global\let\CurrentInput=\INPSTACKB \else%
\ifnum\INPSP=2 \global\let\CurrentInput=\INPSTACKC \else%
\ifnum\INPSP=3 \global\let\CurrentInput=\INPSTACKD \else%
\ifnum\INPSP=4 \global\let\CurrentInput=\INPSTACKE \else%
\ifnum\INPSP=5 \global\let\CurrentInput=\INPSTACKF \else%
               \global\let\CurrentInput=\INPSTACKX \fi\fi\fi\fi\fi\fi%
\WinEdt{<-}%
\global\LASTLINE=\inputlineno%
\global\advance\LASTLINE by -1%
\SRC}
\def\INPUT#1{\relax}
\let\OldINCLUDE=\include
\def\include#1{
\EJECT%
\PUSH{#1.tex}%
\OldINCLUDE{#1}%
\POP}

\let\originalxxxeverypar\everypar
\newtoks\everypar
\originalxxxeverypar{\the\everypar\expandafter\SRC}
\everymath\expandafter{\the\everymath\expandafter\SRC}
\let\zzzxxxbibliography=\bibliography
\def\bibliography#1{\PUSH{\jobname.bbl}\zzzxxxbibliography{#1}\POP}
\output\expandafter{\expandafter\SRCOKfalse\the\output}

\title{INFINITY \vspace{0.2cm} \\ A simple, but not too simple introduction}
\author{by Martin Meyries}

\begin{document}

\maketitle

\subsection*{Introduction}

This text tries to give an elementary introduction to the mathematical properties of infinite sets. The aim is to keep the approach as simple as possible. Advanced knowledge of mathematics is not necessary for a proper understanding, and there is (almost) no use of formulas. At the same time, it is tried to keep the reasoning rigorous and transparent. Major arguments are explained in great detail, while it is explicitly stated if something is not further explained. 

However, also a fully guided trail can become difficult at times. Persistence is necessary to follow. No reader should be discouraged when he or she cannot understand everything on first readings -- infinity is not a piece of cake.
 
The reasoning presented here is mainly due to the famous German mathe\-matician \textsc{Georg Cantor} (1845 -- 1918). Despite the resistance from the mathematical community at that time, he opened the door to infinity almost on his own, dramatically changing the foundations of mathematics. \medskip

The following notions and topics are treated: potential und actual infinity, cardinality and countability, \textsc{Hilbert}'s hotel, \textsc{Cantor}'s diagonal arguments and the highly astonishing answer to the question of the validity of \textsc{Cantor}'s continuum hypothesis.\medskip

\begin{footnotesize} Any kind of questions, comments or suggestions by email is highly welcomed from anyone, especially from non-mathematicians: \href{mailto:m_meyries [at] gmx [dot] de}{m\_meyries $[$at$]$ gmx $[$dot$]$ de}

I am grateful to Alessa Binder, Andreas Bolleyer und Philip M\"uller for their critical reading of an earlier German version (available at \href{http://arxiv.org/abs/1310.4739}{\texttt{arXiv:1310.4739}}), and to Anton Falk for the motivation to produce an English version. 

This article was created on the Hawai'ian island of Kaua'i.\end{footnotesize}

\newpage

\subsection*{Potential infinity}
Sometimes kids claim that there is a largest number. Asking them which one this should be, one receives answers in hundreds, thousands or millions, depending on the knowledge and the age of the kid. It is, however, an insight which is as simple as fundamental that \emph{a largest number does not exist}. In fact, any given number may be increased by 1 to obtain a larger number.

In everyday life one is used to handle numbers up to, say, a million. One is able to imagine concrete things for such numbers: traveling a distance of 100 kilometers, preparing 1500 grams of sugar for a big birthday cake or earning 100.000 Dollars per year. It is already quite difficult to imagine the millions, billions or even trillions being relevant in politics, business, and media.

Even larger numbers, having at least some kind of proper meaning, may often only be found in the natural sciences. For instance, the mass of the earth is about
$6.000.000.000.000.000.000.000.000$  (six septillion) kilograms, while the sun is even one million times heavier. Incredibly large is the total number of elementary particles in the universe, which is nowadays estimated as roughly $10^{87}$ -- this is 1 followed by 87 zeros! But using our imagination, it is simple to overbid even this. We may think of the number $10^{87}+1$ (one more particle), $2\times 10^{87}$ (two universes), or even 1 followed by $10^{87}$ zeros (each particle as a zero). \medskip

We conclude that we can imagine arbitrarily large numbers. In other words, there are \emph{arbitrarily many} natural numbers
$$1,\;2,\;3,\;4,\;5,\;6,\;7,\;8,\;9,\;...\;.$$
One therefore also says that the natural numbers are \emph{potentially infinite}. Their reservoir is inexhaustible and unlimited.\medskip

Potential infinity is in fact a part of reality. Standing on a square and deciding in which direction to go, there are, being pedantic, arbitrarily many possibilities to choose from. Namely, there are as many possibilities as there are angles between 0 and 360 degrees: $1$ degrees, $\nicefrac{1}{2}$ degrees, $\nicefrac{1}{3}$ degrees, $\nicefrac{1}{4}$ degrees, and arbitrarily many degrees more.

\subsection*{Actual infinity}

In this way we do, however, not get in touch with "true infinity", which we will call \emph{actual infinity} in what follows. Despite the fact that the reservoir of natural numbers is unlimited, we can choose, use or imagine only a \emph{finite} number of concrete objects at the same time. Or have you ever really seen, felt or heard infinitely many things \emph{at once}?

The mental collection of well defined and distinct objects to an object in its own right is called a \emph{set}. This naive definition goes back to \textsc{Georg Cantor}. For instance, we may collect the numbers 1, 2 and 3 to form a set, for which we then write $\{1,2,3\}$. In this context, the brackets ``$\{$'' and ``$\}$'' are called set brackets, and each single piece of which a set is formed of is called an \emph{element} of the set. Hence 2 is an element of the set $\{1,2,3\}$.\medskip

Now, actual infinity is something egregious radical. Let us consider the natural numbers as a whole, i.e., \emph{collect all imaginable natural numbers to an object in its own right}. Then the resulting set $\mathbb N$,
$$\mathbb N = \{1,\;2,\;3,\;4,\;5,\;6,\;7,\;8,\;9,\;...\}\,,$$
is not finite, and in this sense it is actual infinite. (Here, the symbol ``$\mathbb N$'' is simply a squiggled ``N'' and is an abbreviation for "natural numbers".) In fact, we call a set \emph{finite}, if we can enumerate all its elements. Concerning $\mathbb N$, we cannot enumerate all its elements -- no matter how far we count, there will always be elements left over. \medskip

Is this acceptable? Are we really \emph{allowed} to collect infinitely many objects to an object in its own right, even if this has no equivalent in reality? This is far from being clear. If at all, actual infinity can only exist in our minds. \medskip

The distinction of potential and actual infinity, which goes back to the famous Greek philosopher \textsc{Aristoteles}, might not be easy to understand. We therefore try to emphasize it once more. On the one hand, for potential infinity one only considers a fixed, finite collection of objects, where there is no upper bound in the number of objects one collects. On the other hand, for actual infinity one directly lifts finiteness and in this sense obtains a non-finite, hence infinite entity.

\subsection*{One small step for \textsc{Cantor}, one giant leap for mankind}
Before \textsc{Cantor}, no mathematician or philosopher seriously dealt with the concept of actual infinity. Even in 1831, one of the greatest mathematicians of all times, \textsc{Carl Friedrich Gau\ss{}}, wrote:

\begin{quote}
"So I heavily protest against the use of an infinite entity as a complete one, which can never be allowed in mathematics."
\end{quote}

By a "complete infinite entity", \textsc{Gau\ss{}} just means the collection of distinct objects to a non-finite whole. An example for this is the set $\mathbb N$ of all natural numbers. \medskip

In the beginning \textsc{Cantor}'s work, starting in 1870, met partially heavy refusal amongst mathematicians. A well known quote, apparently having unknown origin, however, says:

\begin{quote}
"Future generations will consider set theory as a disease, from which one has recovered."
\end{quote}

Today, \textsc{Cantor} is considered as one of the greatest mathematicians ever. \textsc{David Hilbert}, who early recognized the fundamental importance of \textsc{Cantor}'s work, said in 1926:

\begin{quote}
"No one shall ever expel us from the paradise which \textsc{Cantor} has created for us."
\end{quote}

Until today, there are some mathematicians and philosophers who do not accept the concept of actual infinity. In the sequel we do not further want to deal with this question and take the position: why not? Let's see what happens.

\subsection*{Where to begin?}
We thus say "yeah!" to actual infinity. In modern set theory, the existence of infinite sets is in fact formulated as an axiom, i.e., it is considered as fact without further justification. We are going to return to this issue at the end of this text.\medskip

So, how and where do we begin? How can we investigate actual infinity in more detail? How can one investigate an idea or a concept at all? Since we are in uncharted waters and we cannot rely on any kind of experience we can only ask for very basic properties.

We therefore ask: does actual infinity have \emph{any structure at all}? Or is actual infinity just actual infinity, without any further nuances, and this is already the end of the story? In other words:

\begin{quote}
\emph{Are there different types of actual infinity?}
\end{quote}

We shall throughout be guided by this question.

\subsection*{Comparison of infinite sets -- the limits of intuition}
A rather simple observation is that there is not only one infinite set. Besides the already mentioned natural numbers
$$\mathbb N = \{1,\;2,\;3,\;4,\;5,\;6,\;7,\;8,\;9,\;...\}$$  
we may consider the set of all even numbers $\mathbb E$,
$$\mathbb E = \{2,\;4,\;6,\;8,\;10,\;12,\;14,\;16,\;18,\;...\}\,,$$ 
which consists precisely of those natural numbers which do not leave a remainder when divided by 2. We obtain the set $\mathbb E$ if we omit every second number from $\mathbb N$ --  only every second number is even. \medskip

Now, both sets $\mathbb E$ and $\mathbb N$ are of infinite size, and in this sense they are of "equal size". On the other hand, $\mathbb E$ possesses infinitely many elements less than $\mathbb N$, and in this sense it should be "less infinite" than $\mathbb N$. Which point of view is the right one? How can we even \emph{decide} which one is the right point of view? \medskip

We illustrate this problem by means of an example. Imagine that someone possesses an infinite amount of money and has to pay 50$\%$ taxes on it. This surely costs half of the fortune. But after the taxes are paid, there is still an infinite amount of money left. So, is something lost at all? Yes, on the one hand, even an infinite amount! But, on the other hand, nothing is lost, the lucky dude still has an infinity amount of money. Of course, considering a finite fortune, things are totally clear -- paying tax on a million dollar results in possessing less than a million afterwards. \medskip

It seems that one cannot apply this reasoning to infinite sets. Our intuition leads to two contradicting conclusions. We just cannot count the elements of infinite sets and compare the results!

However, not only in mathematics, confusion is often a result of talking about things without being aware of what one is really talking about. Above we used the terms "equal size" and "less infinite" in the context of infinite sets. But what does this mean? When can we say that infinite sets are of "equal size", or that one set is "less infinite" than the other?

\subsection*{Equal cardinality}
How can we compare the "sizes" of infinite sets? \textsc{Cantor}'s solution of this problem is as simple as ingenious. In fact, we do not want to know "how many" elements an infinite set has. This information simply does not exist. We only want to \emph{compare} their sizes, not more and not less.

\textsc{Cantor}'s idea is based on the fact that there is another way to compare the sizes of finite sets which does not ask for the number of elements: \emph{pairing without remainder}.\medskip

Let's consider an illustrating example. Image that we are on the opera ball in Vienna. During a break, men in black tails and women in white ball gowns are mixed up  on the dance floor in a rather confusing way. Now we would like to know whether there are more men or more women on the floor. We could count the women and the men separately, and compare the resulting numbers. However, this gives us much more information than we were originally interested in. In fact, it gives us the precise numbers of women and men.

Another, much simpler possibility is as follows. Using the loudspeakers, we could ask each woman to grab a man. If women are left, then we are sure that there are more women than men on the dance floor. If men are left, there are more men. And if no one is left without a partner, i.e., there is no remainder, then we know that the sets of women and men are of equal size.

This way of comparing sets is in fact well known to kids. They will never portion a box of sweets by counting them and divide the result. Too little kids just cannot do that! And they don't have to. They will rather portion the sweets as follows: one for you, one for me, one for you, one for me, and so on. \medskip

Concerning infinite sets, we are in a similar situation as the little kids: we just can't enumerate the elements. However, we can always try to create pairs of their elements and see if there is a remainder. And this is \textsc{Cantor}'s idea to compare the sizes of infinite sets, resulting in the following \emph{definition}. Instead of the term "size", which might be understood to be based on a finite enumeration, for possibly infinite sets we will rather use the term \emph{cardinality} in everything that follows.

\begin{quote}
Two arbitrary sets $A$ and $B$ are called \emph{of equal cardinality}, if there is at least one way to create pairs from their elements, such that no element from $A$ and $B$ remains left over. If this is the case, we symbolically write $|A| = |B|$.
\end{quote}
In this definition, $A$ and $B$ are symbolic place holders for arbitrary sets of our choice. For instance, we may set $A = \mathbb N$ and $B = \mathbb E$ to compare the cardinality (i.e., the "size") of the set of natural numbers $\mathbb N$ and the set of even numbers $\mathbb E$.\medskip

Before continuing, we first want to reflect on \textsc{Cantor}'s definition.

As we have seen above, for finite sets $A$ and $B$ it leads to nothing but what one usually thinks of when one says that "$A$ and $B$ are of the same size". Moreover, the definition applies to arbitrary sets. We therefore \emph{believe} that the definition makes some sense.

However, at this point it is not clear at all if this really is a "good" definition. It is possible that our notion of equal cardinality doesn't teach us anything on actual infinity. In fact, it might turn out that all infinite sets are of equal cardinality, i.e., we can find a pairing without remainder between all infinite sets. Then the notion would not allow to describe any kind of structures within actual infinity. The reader shall be told in advance that \textsc{Cantor}'s definition is indeed a very good one and yields new insights in actual infinity. \medskip

Finally, we would like to illustrate a pairing without remainder for the finite sets $A = \{1,2,3\}$ and $B= \{a,b,c\}$. This is very simple. Such a pairing is as follows:
$$1\leftrightarrow a,\qquad 2\leftrightarrow b, \qquad 3\leftrightarrow c.$$
But also
$$1\leftrightarrow b,\qquad 2\leftrightarrow c, \qquad 3\leftrightarrow a,$$
describes a pairing without remainder. Clearly, in this case it is lot quicker to realize that $A$ and $B$ both have three elements and are thus of equal cardinality or size. But for infinite sets we cannot do that, the only thing we can rely on is our definition.

\subsection*{The equal cardinality of $\mathbb N$ und $\mathbb E$}
We can indeed find a pairing of the elements of the natural numbers $\mathbb N$ and the even numbers $\mathbb E$, without leaving a remainder. This is based on the following idea of \emph{infinitely counting} the even numbers: we may consider 2 as the \emph{first} even number, 4 as the \emph{second} even number, 6 as the \emph{third} even number, and so on. In a scheme, the pairing looks as follows:

$$\begin{array}{rcccccccccc}
\mathbb N :\;\;& 1 & 2 & 3 & 4 & 5 & 6 & 7 & 8 &9  &\cdots\\
& \updownarrow & \updownarrow &\updownarrow &\updownarrow &\updownarrow &\updownarrow &\updownarrow &\updownarrow &\updownarrow &\\
\mathbb E :\;\;& 2 & 4 & 6 & 8 & 10 & 12 & 14 & 16 & 18  &\cdots
\end{array}
$$
In other words: in our pairing, the partner of an element of $\mathbb N$ is its double value -- which surely is an even number and thus an element of $\mathbb E$. In this way, no element of $\mathbb E$ is left over, the partner of a given even number is just its half value.

Following \textsc{Cantor}'s definition, we conclude that $\mathbb N$ und $\mathbb E$ are of equal cardinality,
$$|\mathbb N| = |\mathbb E|.$$
If we believe that the definition makes sense, as indicated above, then we must accept that there are as many natural numbers as even numbers, despite the fact that $\mathbb E$ is a strict subset of $\mathbb N$. We may further conclude the following:

\begin{quote}
\emph{By removing infinitely many elements from an infinite set, its cardinality can remain equal.}
\end{quote}

The latter is impossible in the case of finite sets: by removing elements from a finite set, the cardinality does not remain equal, as the interested reader may try to prove on his or her own. We see that in dealing with infinite sets we immediately get in touch with new and unknown phenomena. We remove every second element of $\mathbb N$, and still do not change its cardinality (its "size"). 

Of course, this does not mean that \emph{no matter} how many elements we take from an infinite set, the cardinality always remains equal. If we remove every element from $\mathbb N$ except the element 1, then the sets $\{1\}$ and $\mathbb N$ are not of equal cardinality (why?). But what if we leave infinitely many elements left over? We will soon return to this question.\medskip

First, another important remark is in order. A further new phenomenon concerning infinite sets is that even though we have successfully found a pairing without remainder, there might be other pairings which leave a remainder, even one which is infinitely large. For instance, the pairing
$$\begin{array}{rcccccccccc}
\mathbb N :\;\;& 1 & 2 & 3 & 4 & 5 & 6 & 7 & 8 &9  &\cdots\\
&  & \updownarrow & &\updownarrow & &\updownarrow & &\updownarrow & &\\
\mathbb E :\;\;&  & 2 &  & 4 &  & 6 & & 8 &   &\cdots
\end{array}
$$
misses all odd numbers of $\mathbb N$. However, this phenomenon does not contradict our definition of equal cardinality. It only asks for \emph{at least one} pairing without remainder. It doesn't matter whether there are other pairings with remainder.

\subsection*{Larger and smaller cardinalities}

Based on the same idea as before, we may define when an arbitrary set is of larger cardinality (or simply larger) than another set.

\begin{quote}
The set $B$ is called \emph{of larger cardinality} than the set $A$, if $A$ and $B$ are not of equal cardinality and if $A$ is of equal cardinality as a subset of $B$. In this case one symbolically writes $|B| > |A|$.
\end{quote}
Of course, for $|B| > |A|$ one may equivalently write $|A| < |B|$ and say that $A$ is of \emph{smaller cardinality} than $B$.\medskip

Let us pause for a second and have a closer look at this definition. To call $B$ of larger cardinality than $A$, two properties must be satisfied. First, $A$ and $B$ should not be of equal cardinality, i.e., every pairing of elements leaves a remainder. It is clear that this should be satisfied, otherwise it would be rather strange to call $B$ larger than $A$. The second property is designed to decide which of the sets is the larger one: if we remove a suitable part of $B$, i.e., we only consider a certain subset of $B$, then we should obtain a set of equal cardinality as $A$. In this case it is fair enough to call $B$ larger than $A$. \medskip

Let us illustrate all this in the seemingly simple case of the finite sets $A = \{1,2,3\}$ and $B = \{a, b, c, d\}$. We immediately see that $B$ is larger than $A$ in the usual sense, since it contains one element more than $A$ does. But as before, this not the point here. We only want to rely on our definitions, as this is the only firm ground we can rely on when it comes to dealing with infinite sets! 

We want to prove in the sense of \textsc{Cantor} that $B$ is larger than $A$. We have seen before that $A$ is of equal cardinality as the subset $\{a,b,c\}$ of $B$. Hence we are done with the second required property. Now we only have to prove that $A$ and $B$ are not of equal cardinality, i.e., no pairing without remainder is possible. This is clear, isn't it? Well, as mathematicians we fully want to understand everything we are talking about and we want to be 100$\%$ sure that all our claims are true. Second, nothing will be clear anymore when dealing with infinite sets. It is thus a good idea to see how one can exclude the possibility of a pairing without remainder already in a simple case.

A foolproof method is to make the effort to go through all possible pairings, which are as many as $24 = 4\times 3\times 2$, and to check whether there always is a remainder. An example is
$$1\leftrightarrow c, \qquad 2\leftrightarrow d,\qquad 3\leftrightarrow a, \qquad ?\leftrightarrow b.$$
Here the element $b$ of $B$ remains without a partner from $A$. Now the reader may convince his- or herself that the same happens for all of the other 23 pairings. We therefore conclude from our definition alone, that $B$ is of larger cardinality than $A$.

This example shows that it may be a big effort to check whether a set is larger than another one. However, in the case of finite sets one does not have to check every possible pairing as it is indicated above. One can show (we will not do so here), that all pairings lead to the same result -- either all pairings are without a remainder, or all pairings leave a remainder, and the latter is always of the same size. But to emphasize it once more: for infinite sets there might not be such a shortcut such that we can only use the foolproof method. \medskip

One should moreover emphasize again that at this point we only have a \emph{method} to detect different types of actual infinity. It is not clear yet whether there really exist infinite sets with non-equal cardinality, or if all infinite sets are of the same cardinality.

\subsection*{The cardinality of $\mathbb N$ is the smallest infinite cardinality}

We have seen above that even if we remove every second number from the set of natural numbers $\mathbb N$, the resulting set $\mathbb E$ of even numbers is of equal cardinality as $\mathbb N$. Now we want to go further and ask what happens if we remove even more numbers from $\mathbb N$\,? Do we eventually obtain an infinite set, which is of smaller cardinality than $\mathbb N$\,? In full generality, we ask: does there exist an infinite set which is of smaller cardinality than the set of natural numbers?\medskip

Here the answer is no. To see that, \emph{assume} that there is an arbitrary infinite set $A$ which is \emph{not} of equal cardinality as $\mathbb N$. (This does not imply that such a set $A$ really exists, we are only assuming it!) We have to keep in mind that we really know nothing else about $A$, except that it is an infinite set and not of equal cardinality as $\mathbb N$. 

The first information is enough to create the following pairing: as the partner of $1$ we choose an arbitrary element of $A$, as the partner of 2 we choose an arbitrary further element of $A$, as the partner of 3 we choose another arbitrary further element of $A$, and so on. In this way, each element of $\mathbb N$ obtains a partner from $A$. Why is that the case? If, for instance, from the number 1000 on there would be no further elements left from $A$, then $A$ would have only had 1000 elements, which contradicts our assumption that it has infinitely many elements.

In this way we obtain a pairing without remainder between $\mathbb N$ and the set of elements of $A$, which were chosen in this process. Hence $\mathbb N$ and this subset of $A$ are of equal cardinality, as a pairing without remainder is precisely what our definition asks for. Altogether, $\mathbb N$ and $A$ are not of equal cardinality (this was our second assumption above) and $\mathbb N$ is of equal cardinality as a subset of $A$. Hence, according to the above definition, the set $A$ (if it really exists) is of larger cardinality than $\mathbb N$. We conclude:

\begin{quote}
\emph{Every infinite set is of equal or larger cardinality as $\mathbb N$.}
\end{quote}

So, \emph{if} there are different infinite cardinalites, then there is a smallest one, namely the one of the set of natural numbers. This is a truly fundamental finding. The cardinality of $\mathbb N$ is something very special. Hence one would like to term it separately.

\subsection*{Countability}

As mentioned before, one may interpret the pairing without remainder between $\mathbb E$ and $\mathbb N$ as an "infinite counting" of the elements of $\mathbb E$. One considers the partner of 1 as the \emph{first} even number, the partner of 2 as the \emph{second} even number, the partner of 3 as the \emph{third} even number, and so on. This also means that the elements of any set which is of equal cardinality as $\mathbb N$, can be written in an infinite list: the partner of 1 as the first entry, the partner of 2 as the second entry, and so on.\medskip

One therefore calls a set of equal cardinality as the natural numbers \emph{infinitely countable} or \emph{infinitely listable}. For instance, the set $\mathbb E$ is infinitely countable, and soon we will find more infinitely countable sets.

Often one omits the word "infinitely" and simply says \emph{countable} instead of infinitely countable. The reader should not be confused by that. Of course, it is impossible to count or enumerate infinitely many objects in the usual sense. This is only possible for finite sets. If at all, then only an infinite counting in the sense of a pairing with $\mathbb N$ without remainder is possible for an infinite set.\medskip

A set, which is of larger cardinality than $\mathbb N$, is not infinitely countable or listable. We therefore want to call such a set \emph{infinitely uncountable} or just \emph{uncountable}. Note that, once more, at this point we do not know whether an uncountable set exists.

\subsection*{Further countable sets}

We have seen that the cardinality of the set $\mathbb N$ of natural numbers is the smallest infinite cardinality. Now we want to head to another direction and add elements to $\mathbb N$, hoping to eventually find a set of larger cardinality, i.e., an uncountable set.

Taking a back seat, we start with only a single element and consider

$$\mathbb N_0 = \{0,\;1,\;2,\;3,\;4,\;5,\;6,\;7,\;8,\;...\},$$
which is nothing but the set of natural numbers $\mathbb N$ together with zero.  A pairing without remainder can be realized as follows:

$$\begin{array}{rcccccccccc}
\mathbb N :\;\;& 1 & 2 & 3 & 4 & 5 & 6 & 7 & 8 &9  &\cdots\\
& \updownarrow & \updownarrow &\updownarrow &\updownarrow &\updownarrow &\updownarrow &\updownarrow &\updownarrow &\updownarrow &\\
\mathbb N_0 :\;\;& 0 & 1 & 2 & 3 & 4 & 5 & 6 & 7 & 8  &\cdots
\end{array}
$$
In short, the pairing rule says that one simply has to subtract 1 from an element of $\mathbb N$ to find its partner in $\mathbb N_0$. We conclude that $\mathbb N$ and $\mathbb N_0$ are of equal cardinality, for which we symbolically write $|\mathbb N| = |\mathbb N_0|$. Hence, also $\mathbb N_0$ is a countable set.

In the same way finds that adding 19 or a million new elements does not enlarge the cardinality of $\mathbb N$. Here the rule for pairing without remainder is: subtract 19 or one million to find the partner. We conclude that by adding only a finite number of elements the cardinality of $\mathbb N$ does not change. \medskip

In a next step we add infinitely many elements to $\mathbb N$ and consider the set $\mathbb Z$ of all \emph{integers},
$$\mathbb Z = \{..., \;-4,\;-3,\;-2,\;-1,\;0,\;1,\;2,\;3,\;4,\;...\}.$$
In addition to all natural numbers and zero, the integers contain all negative numbers as well. Now we might have found an uncountable set, since we see that $\mathbb Z$ is "unlimited to the left and to the right", while $\mathbb N_0$ and $\mathbb E$ are only "unlimited into one direction". But it turns out that this is not the point. Writing $\mathbb Z$ in a different way, namely
$$\mathbb Z = \{0,\;1,\;-1,\;2,\;-2,\;3,\;-3,\;4,\;-4,\;...\},$$
it becomes clear how a pairing with $\mathbb N$ without remainder works out:
$$\begin{array}{rcccccccccc}
\mathbb N :\;\;& 1 & 2 & 3 & 4 & 5 & 6 & 7 & 8 &9  &\cdots\\
& \updownarrow & \updownarrow &\updownarrow &\updownarrow &\updownarrow &\updownarrow &\updownarrow &\updownarrow &\updownarrow & \\
\mathbb Z :\;\;& 0 & 1 & -1 & 2 & -2 & 3 & -3 & 4 & -4  &\cdots
\end{array}
$$
Let us consider the pairing rule in more detail, as it is slightly more involved than our previous pairings. It distinguishes between even and odd elements of $\mathbb N$. An even number is divided in half, for instance, the partner of the element 4 of $\mathbb N$ is the element 2 of $\mathbb Z$. For an odd element of $\mathbb N$, one subtracts 1, divides in half and finally puts a minus in front of it. For instance, this yields $-4$ as the partner of $9$, since $-4 = -(9-1):2$. To be sure that we really catch every element of $\mathbb Z$, we invert these rules to find the corresponding partner in $\mathbb N$. A positive number is doubled; for a negative number we omit the minus, then double it and finally add 1.\medskip

We may conclude that $\mathbb Z$ is countable as well. There are as many integers as natural numbers! This was somehow to expect after we had rewritten $\mathbb Z$. The step from $\mathbb Z$ to $\mathbb N$ is analogous to the step from $\mathbb N$ to $\mathbb E$, namely removing every second element, and we had already seen that $\mathbb N$ and $\mathbb E$ are of equal cardinality.

We record the following fundamental general consequence:

\begin{quote}
\emph{The cardinality of an infinite set is not necessarily enlarged when adding infinitely many new elements.}
\end{quote}

The fact that all  infinite sets we have investigated so far are countable can symbolically be written as follows:

$$|\mathbb N| = |\mathbb E| = |\mathbb N_0| = |\mathbb Z|.$$

\subsection*{An excursion to \textsc{Hilbert}'s hotel}

Before continuing in investigating larger and larger sets, aiming to find a larger cardinality than the one of $\mathbb N$, we would like to present a very interesting \emph{Gedankenexperiment} (thought experiment) which is due to \textsc{David Hilbert} and nowadays known as \textsc{Hilbert}'s hotel. It illustrates some possibly strange consequences that arise when one mixes up everyday life with the properties of actual infinite sets. \medskip

Imagine there is a hotel with infinitely countable many rooms, i.e., as many rooms as there are natural numbers. Further imagine that all these rooms are occupied and there is new guest asking for vacancy. Do we have to send the guest to another hotel because we cannot offer a room? At first sight, this is the case, since all rooms are occupied. \medskip

However, \textsc{Hilbert}'s hotel is not an ordinary, finite hotel. In the previous section we realized that the cardinality of $\mathbb N$ does not change if we add a single element. In our pairing rule every natural number simply moved by 1 to the right and after that there was space for another element, namely zero.

In the same way we can create vacancy for the new guest in our infinite hotel: we ask the guest in room 1 to move to room 2, the guest from room 2 to move to room 3, the guest from room 3 to move to room 4, and so on. Then all previous guests are again in a room -- and room 1 is vacant for the new guest!

Of course, we can proceed in a similar way if 346 or one billion new guests arrive and ask for vacancy.\medskip

So, what happens if \emph{infinitely countable many new guests} ask for vacancy? Even in this case we don't have to send anyone to another hotel. Inspired by the pairing rule between $\mathbb N$ and $\mathbb E$, we ask all of the present guests to move into the room with the double room number. After that, only the rooms with even room numbers are occupied. Hence all rooms with odd room number are vacant and we have enough space for as many guests as there are odd numbers. Since there are as many odd numbers as natural numbers (the reader should try to find the pairing rule in this case!) there is enough space for countable many new guests.\medskip

To summarize: a hotel with (countable) infinitely many rooms always has vacancy, even for (countable) infinitely many new guests.

Of course, this does not match our "finite" intuition. But that's no problem: in reality, there just are no hotels with infinitely many rooms. We should not forget from where we started, namely the assumption (the axiom) that there are actual infinite sets. This assumption only exists as a concept or idea and has no counterpart in reality. It is therefore not surprising that this assumption has consequences which only exist in the world of ideas as well.\medskip

Finally, some guests might be annoyed because they are asked to move so often and thus give bad online ratings for \textsc{Hilbert}'s hotel. However, we can offer all rooms even for free. How is this possible? If we need some money, we ask the guest in room 1 for it, who asks the guest in room 2 for the same amount, and so on...

\subsection*{Countability of the rational numbers -- \\\textsc{Cantor}'s first diagonal argument}
We continue our search for an uncountable set. Together with all negative integers we further add all fractions to $\mathbb N$, which yields the set of \emph{rational numbers} $\mathbb Q$,
$$\mathbb Q = \left \{\nicefrac{p}{q}\;:\; \text{$p$ and $q$ from $\mathbb Z$ with $q$ not equal to $0$}\right\}.$$
For instance, $\nicefrac{1}{3}$, $-\nicefrac{7}{28}$ and $\nicefrac{1000}{9}$ are rational numbers. All integers are rational as well, since, for instance, we have $5 = \nicefrac{5}{1}$ and $-2 = \nicefrac{-2}{1}$.

By the way, here one uses a squiggled "Q" as an abbreviation for the set of rational numbers since the squiggled "R", namely $\mathbb R$, is already reserved for the even more important set of real numbers, which will soon enter the scenery.\medskip

Now we finally seem to have found a new type of actual infinity. There are "gaps" between the elements of the sets that we have considered so far: between 4 and 5 there is no further natural number. This is entirely different for the rationals: between two arbitrary rational numbers there are even infinitely many more rational numbers! One can understand this by means of an example: between 0 and 1 there is $\nicefrac{1}{2}$, between 0 and $\nicefrac{1}{2}$ there is $\nicefrac{1}{3}$, between 0 and $\nicefrac{1}{3}$ there is $\nicefrac{1}{4}$, and so on. Hence $\mathbb Q$ is rather densely packed. It appears that there are a lot "more" fractions than only natural numbers.\medskip

But also here we are mislead by our intuition. \textsc{Cantor} has proved that the set of rational numbers $\mathbb Q$ is countable as well.\medskip

What does \textsc{Cantor}'s pairing rule without remainder look like? For simplicity we only want to consider a pairing between $\mathbb N$ and the positive rational numbers, and leave it to the reader to how to conclude from this that $|\mathbb Q| = |\mathbb N|$.

\textsc{Cantor}'s ingenious idea how to infinitely count the rational numbers is based on the following infinite rectangular scheme,

$$
\begin{array}{cccccc}
\nicefrac11  &   \nicefrac12  & \nicefrac13  &\nicefrac14      & \cdots \\\\
\nicefrac21  & \nicefrac22  & \nicefrac23  &\nicefrac24      & \cdots \\\\
\nicefrac31  & \nicefrac32  & \nicefrac33  &\nicefrac34     & \cdots \\\\
\nicefrac41  & \nicefrac42  & \nicefrac43  &\nicefrac44 & \cdots \\
\vdots  & \vdots  & \vdots   & \vdots &  &
\end{array}
$$
In this scheme, every positive rational numbers appears at least once. Of course, we do not obtain a pairing between $\mathbb N$ and $\mathbb Q$ without remainder if we start in the upper left corner at $\nicefrac11$ and simply count into a vertical or horizontal direction. In this way we wouldn't even come close to most of the fractions. It is a much better idea to count \emph{diagonal} --  which precisely is \textsc{Cantor}'s \emph{first diagonal argument}. It works as follows:
$$
\begin{array}{ccccccccc}
\nicefrac11  &\rightarrow &   \nicefrac12  & &\nicefrac13  & \rightarrow & \nicefrac14      & \cdots \\
& \swarrow & & \nearrow & & \swarrow & & \iddots\\
\nicefrac21  & &\big (\nicefrac22\big)  & &\nicefrac23  & &\big (\nicefrac24   \big)   & \cdots \\
\downarrow & \nearrow & & \swarrow & & \nearrow & & \iddots \\
\nicefrac31  & &\nicefrac32  & &\big(\nicefrac33\big)  & &\nicefrac34     & \cdots \\
& \swarrow & & \nearrow  &&\nearrow && \iddots\\
\nicefrac41  & &\big (\nicefrac42 \big) &  &\nicefrac43  && \big(\nicefrac44\big) & \cdots \\
\vdots  &\iddots & \vdots  & \iddots&\vdots   & \iddots & \vdots &  &
\end{array}
$$
Let us explain it in more detail. We start in the left upper corner, the partner of 1 is $\nicefrac11 = 1$. The partners of 2 and 3 are the two entries on the second diagonal, namely $\nicefrac12$ and $\nicefrac21$, respectively. The partners of 4 and 5 are the entries $\nicefrac31$ and $\nicefrac13$ on the third diagonal, respectively. Note here that $\nicefrac22 = 1$ does not need a partner in $\mathbb N$, since we have already defined a partner for 1. Next, the partners of 6, 7, 8 and 9 are the entries $\nicefrac14$, $\nicefrac23$, $\nicefrac32$ and $\nicefrac41$ on the fourth diagonal, respectively. And so on. In this way we are snaking along the diagonals through the whole scheme and find for each element of (the positive subset of) $\mathbb Q$ a partner from $\mathbb N$, and vice versa.

We conclude that the rational and the natural numbers are of equal cardinality,
$$|\mathbb N| = |\mathbb Q|.$$

\subsection*{The continuum -- the set of real numbers}
Now we stop pulling punches and consider \emph{all decimal numbers at once}. Examples are all rational numbers, for instance
$$\nicefrac13 = 0,33333...\;, \qquad 5 = 5,00000...\;, \qquad -\nicefrac65 = -1,20000...\;;$$
further all roots of all rational numbers, e.g.,
$$\sqrt{2} = 1,41421...\;, \qquad -\sqrt[3]{\nicefrac{1001}5} = -5,84998...\;$$
and finally all remaining decimal numbers, the so called \emph{transcendental numbers}, e.g., fabulous $\pi$ and \textsc{Euler}'s number $e$,
$$\pi = 3,14159...\;,\qquad e = 2,71828...\;.$$
The resulting set is called the set of \emph{real numbers} and one writes $\mathbb R$ for this set. Similar as for $\mathbb N$, $\mathbb Z$ and $\mathbb Q$, one can find an expression for $\mathbb R$ in terms of set brackets, $\mathbb R = \{ ... \}$. But this looks rather complicated and is not so important for what follows. What one should have in mind is simply that $\mathbb R$ consists precisely of all decimal numbers.\medskip

In particular, $\mathbb Q$ is a subset of $\mathbb R$, as we added even infinitely many elements to obtain $\mathbb R$ from $\mathbb Q$. But what is \emph{really happening} when enlarging $\mathbb Q$ to $\mathbb R$\,? This cannot be described in an elementary way and is not easy to understand. 

As we have seen, $\mathbb Q$ is rather densely packed. However, still $\mathbb Q$ is laced with infinitely many, \emph{infinitely small} "holes". The mathematicians also say that $\mathbb Q$ is not complete. For instance, at $\sqrt{2}$ there is a "hole" in $\mathbb Q$, in the following sense. One can approximate $\sqrt 2= 1,41421...$ arbitrarily close by rational numbers, e.g., by the sequence
$$1\;;\qquad 1,4\;; \qquad 1,41\;;\qquad 1,414\;;\qquad 1,4142\;;\qquad \text{etc.}$$
But since $\sqrt 2= 1,41421...$ itself is not a rational number  one jumps \emph{out of $\mathbb Q$} in infinitely many steps, observing the hole $\sqrt 2$ at the place where one is falling out.

Only by taking all decimal numbers together we close all of these holes, since one cannot jump out of $\mathbb R$ as described above (each such process ends up on a real number again -- the reader should think about this in more detail). Therefore, since $\mathbb R$ has no holes and one can slide smoothly on it, one calls the set of real numbers the \emph{continuum}.

\subsection*{The real numbers are uncountable -- \\ Cantors second diagonal argument}
As we have seen, in enlarging $\mathbb Q$ to $\mathbb R$ something is happening at an \emph{infinitely small scale}; and this has tremendous consequences at the infinitely large scale. With a further stroke of genius, \textsc{Cantor}
could show that the set of real numbers is uncountable, i.e., that $\mathbb R$ is of larger cardinality than $\mathbb N$. Before considering the beautiful argument, his \emph{second diagonal argument}, in more detail, we recall our guiding question from the beginning of this text, whether actual infinity has any kind of structure. Now we can answer this question with a "yo!":

\begin{quote}
\emph{There are different infinite cardinalities.}
\end{quote}
The significance of \textsc{Cantor}'s discovery cannot be overestimated. According to \textsc{Bertrand Russell} in 1902,
\begin{quote}
"The solution of the difficulties, which earlier encompassed the mathematical infinite, is probably the greatest achievement of which our age can boast."
\end{quote}

How could \textsc{Cantor} prove that $\mathbb R$ is of larger cardinality than $\mathbb N$\,? Recall that two things have to be verified: first, that $\mathbb N$ is of equal cardinality as a subset of $\mathbb R$, and, second, that $\mathbb N$ and $\mathbb R$ itself are not of equal cardinality. The first point is easy to verify, since $\mathbb N$ itself is a subset of $\mathbb R$ (recall here that $2 = 2,00000...$ and so on).

The difficult issue is thus to prove that $\mathbb N$ and $\mathbb R$ are not of equal cardinality. Now we can only rely on the foolproof method and show that \emph{any} pairing of $\mathbb N$ and $\mathbb R$ leaves \emph{at least} one element as a remainder.

How can one prove that? In our previous example with the finite sets $A = \{1,2,3\}$ and $B=\{a,b,c,d\}$ one only has to check 24 pairings, but in case of $\mathbb N$ and $\mathbb R$ there are infinitely many pairings!

A strategy to circumvent this problem is as follows. We are going to pick an \emph{arbitrary} such pairing and \emph{concretely point out} an element of $\mathbb R$ which does not have a partner in $\mathbb N$. If we can achieve this, then we can be sure that \emph{no} pairing without remainder exists. \medskip

\textsc{Cantor}'s idea to find for a given pairing between $\mathbb N$ and $\mathbb R$ a decimal number which does not have a partner from $\mathbb N$, is as follows. Pick an arbitrary such pairing, e.g.,
$$
\begin{array}{ccc}
1 & \leftrightarrow & 0,3333...\\
2 & \leftrightarrow & 0,5432...\\
3 & \leftrightarrow & 0,6775...\\
4 & \leftrightarrow & 0,1010...\\
\vdots & \vdots & \vdots
\end{array}
$$
We construct the desired decimal number without partner in jumping along the diagonal on the right hand side of this scheme, adding 1 to each diagonal element. 

Considering the partner of 1, which is $0,3333...$, we pick the \emph{first} decimal place and add 1, obtaining the \emph{first} decimal place of our desired decimal number, i.e., $3+1=4$. From the partner of 2, which is $0,5432...$, we add 1 to the \emph{second} decimal place, obtaining $4+1=5$, and this becomes the \emph{second} decimal place of the number we are constructing. Its \emph{third} decimal place becomes the \emph{third} decimal place plus 1 of the partner $0,6775...$ of 3, i.e., $7+1=8$. And so on (if a 9 occurs, we make a 0 out of it). This algorithm generates a decimal number which cannot have a partner in the above given pairing. Let us illustrate the construction of this number in our scheme:
$$
\begin{array}{ccc}
1 & \leftrightarrow & 0,\,\framebox[1.05\width]{\textbf{3}}\,333...\\
2 & \leftrightarrow & 0,5\,\framebox[1.05\width]{\textbf{4}}\,32...\\
3 & \leftrightarrow & 0,67\,\framebox[1.05\width]{\textbf{7}}\,5...\\
4 & \leftrightarrow & 0,101\,\framebox[1.05\width]{\textbf{0}}\,...\\
\vdots & \vdots& \vdots\\
? & \leftrightarrow & 0, \textbf{4581...}\\
\end{array}
$$
The resulting decimal number is $0,4581...$\,. Now, can 1 be the partner of this number in the pairing? No, since $0,4581...$ does not coincide with $0,3333...$ in at least its \emph{first} decimal place. Can 2 be its partner? Surely not, since it differs from $0,5432...$ at least in the \emph{second} decimal place. And so on, the number $0,4581...$ is precisely made in a way that it differs from each decimal number in at least one decimal place. Hence the constructed number does not appear in the above list and thus has no partner within the pairing.\medskip

We conclude that in an arbitrary pairing of $\mathbb N$ and $\mathbb R$, at least one element of $\mathbb R$ does not have a partner. A pairing without remainder is therefore impossible and $\mathbb R$ is of larger cardinality than $\mathbb N$. In other words: \emph{The set of real numbers is uncountable.}\medskip

By the way, \textsc{Cantor}'s argument even shows a little more. As we only considered decimal numbers starting with 0, already the set of all real numbers between 0 and 1 is uncountable. More generally, the set of real numbers between two arbitrary numbers is already uncountable. (The reader should convince his- or herself of this proposition.) This shows from another perspective that the crucial thing in the transition from $\mathbb Q$ to $\mathbb R$ happens in the infinitely small.

\subsection*{Infinitely many infinite cardinalities}
We have seen that there are at least two different infinite cardinalities, namely the one of $\mathbb N$ and the one of $\mathbb R$. Without going into the details (even though this issue is not too difficult) we would like to mention that there is even a whole \emph{sequence} of different infinite cardinalities which become larger step by step. We may thus conclude:

\begin{quote}
\emph{There are infinitely many infinite cardinalities.}
\end{quote}

\subsection*{A glimpse on the fine structure of actual infinity -- \\ Cantor's continuum hypothesis}

In this final section we would like to present two further issues concerning the fine structure of actual infinity, avoiding the rather difficult details. Here we are going to reach fundamental limits of knowledge and logic!\medskip

First the natural question arises, \emph{how many} different infinite cardinalities do exist? Are there countable many, as many as real numbers, or even more?

The answer is intriguing: the entity of all different infinite cardinalities is so extremely large that one cannot collect them all at once to a set without generating disastrous contradictions. Already when we collected all natural numbers to a set we saw that we should at least be very careful -- even though today the vast majority of mathematicians believes that this is not problematic. The collection of different objects to a set is however limited, and the entity of all infinite cardinalities is far beyond the limit. \medskip

Things become even more interesting in the vicinity of the cardinality of the set of natural numbers $\mathbb N$. We already know that this is the smallest infinite cardinality. We further know that the set $\mathbb R$ of real numbers, the continuum, is of larger cardinality than $\mathbb N$. It is therefore natural to ask:

\begin{quote}
\emph{Is there a cardinality \emph{in between} the ones of $\mathbb N$ and $\mathbb R$\,?}
\end{quote}
This is to ask whether any infinite subset of $\mathbb R$ is either of the same cardinality as $\mathbb N$ or of the same cardinality as $\mathbb R$.\medskip

\textsc{Cantor} was the first to ask this fundamental question. He, however, couldn't answer it neither in the positive nor in the negative. But he was quite convinced that there is \emph{no} such cardinality. In 1878 he published this conjecture, which is nowadays known as the \emph{continuum hypothesis}.

For many decades, the answer to the above question was unknown, mathematicians and philosophers were not able either to prove or to disprove the continuum hypothesis. In 1900, during the International Congress of Mathematicians in Paris, \textsc{Hilbert} presented a list containing the at that time 23 most important unsolved mathematical problems (this is comparable to the seven \emph{millennium problems} from 2000). In this list, the continuum hypothesis was at first place, showing its utmost significance.\medskip

The final answer, found by \textsc{Kurt G\"odel} (1938) and \textsc{Paul Cohen} (1960), is sensational.

Let us have a closer look at the above question. It asks whether \textsc{Cantor}'s hypothesis is \emph{true or false}. In asking this, one therefore assumes that the hypothesis \emph{really is either true of false}. The reader might cry, "come on, what else should it be?" The assumption that every hypothesis is either true or false is a fundamental instrument for a mathematician, comparable to the telescope for astronomers or the test tube for chemists. However, the following holds true:

\begin{quote}
\emph{The continuum hypothesis is neither provable nor disprovable, neither true nor false, it is \emph{undecidable}.}
\end{quote}
This was definitely not to expect. How is it possible that a hypothesis (an assertion or a proposition) is neither true nor false?

To understand this issue one has to realize what it really means that a mathematical assertion is provable or disprovable. In fact, \emph{proving} an assertion means to logically deduce it from assumptions which are not further justified -- the axioms. On the other hand, to \emph{disprove} an assertion means to logically deduce its contrary from the axioms.

We see that whether an assertion is decidable (i.e., it is either provable or disprovable) depends on two things: the axioms and the logic deduction rules. It is however not clear whether in a system of axioms \emph{any assertion} is really decidable. Until 1931, the mathematicians tacitly assumed this as a given fact.\medskip

In our case, the underlying axioms are the \emph{axioms of set theory}. We have already encountered one of these axioms, namely the \emph{axiom of infinity}, claiming that there exists a set which contains all natural numbers. Further axioms of set theory are, for instance, that there exists a set which does not contain any elements (the \emph{empty set}) and that, given two sets, one can always create their union to obtain a new set. It is moreover a consequence of the \emph{power set axiom} that one may collect all real numbers at once in a set. We have also encountered what the axioms of set theory do not allow, among other things: to collect all different infinite cardinalities to a set.\medskip

The undecidability of the continuum hypothesis must now be understood as follows. One can \emph{prove} (\textsc{G\"odel}) that it is impossible to deduce the contrary of the hypothesis from the axioms of set theory. But one can also \emph{prove} (\textsc{Cohen}) that one cannot deduce the hypothesis itself from the axioms of set theory. 

This means that the continuum hypothesis is \emph{independent} of the axioms of set theory. In principle, now one has the choice whether one adds it as \emph{a new axiom} to the other axioms, whether one adds its negation as a new axiom -- or whether one simply ignores it. \medskip

The undecidability of \textsc{Cantor}'s hypothesis was the first relevant example of an even more fundamental problem in mathematics. Indeed, \textsc{G\"odel} showed in 1931 in his \emph{Unvollst\"andigkeitssatz}, among other things, that most of the systems of axioms which are relevant for mathematics contain undecidable assertions. More precisely, this applies to systems of axioms which contain sufficient information to describe the usual addition and multiplication of natural numbers. These aren't really strong requirements. \textsc{G\"odel}'s results shattered the foundations of mathematics and logic.\medskip

Until the undecidability of the continuum hypothesis was established in 1960, only "artificial" undecidable assertions were known. (An example is the assertion "This assertion is false". Indeed, supposing that the assertion "This assertion is false" is true leads to a contradiction to the assertion itself, which claims the assertion is false. Supposing, on the other hand, the assertion "This assertion is false" to be false also leads to a contradiction, since the contrary of this assertion just means that the assertion "This assertion is false" is true. Hence this assertion is neither provable nor disprovable. Still there?) \medskip

The nowadays widely accepted system of axioms of set theory, as well as the rules of deduction, are man made and should thus not be considered as forever carved in stone. Until today a possible modification of the classical system of axioms of set theory -- possibly deciding the continuum hypothesis -- is an active field in mathematical research. It is very exciting how this story continues. For those who want to know more or even take an active part in it: study mathematics at a university!

\subsection*{Some further literature}
\begin{itemize}
\item K. Devlin, \emph{The Joy of Sets (2nd ed.).} Springer Verlag, 1993.
\item J. Ferreir\'os, \emph{Labyrinth of Thought: A history of set theory and its role in modern mathematics.} Birkh\"auser, 2007.
\item D. R. Hofstadter, \emph{G\"odel, Escher, Bach: An eternal golden braid}. Basic Books, 1999.
\item P. Johnson, \emph{A History of Set Theory.} Prindle, Weber \& Schmidt, 1972.
\item K. Kunen, \emph{Set Theory: An Introduction to Independence Proofs.} North Holland, 1980.
\item M. Potter, \emph{Set Theory and Its Philosophy: A Critical Introduction.} Oxford University Press, 2004.
\item M. Tiles, \emph{The Philosophy of Set Theory: An Historical Introduction to Cantor's Paradise.} Dover Publications, 2004.
\item O. Deiser,  \emph{Einf\"uhrung in die Mengenlehre: Die Mengenlehre Georg Cantors und ihre Axiomatisierung durch Ernst Zermelo.} Springer-Verlag, 2010 (in German).
\end{itemize}

\end{document}